\documentstyle [11pt,amstex]{article}

\newtheorem{thm}{Theorem}

\newtheorem{lemma}{Lemma}
\newtheorem{cor}{Corollary}
\newtheorem{main}{Main Theorem}

\newcommand{\mathbox}[1]{\makebox(0,0){${\scriptstyle #1}$}}
\newcommand{\mathbbox}[1]{\makebox(0,0){${\textstyle #1}$}}
\newcommand{\mtiny}{\scriptscriptstyle}

\newcommand{\proof}{{\it Proof: }}
\newcommand{\qed}{\hfill $\Box$ \vspace{1cm}}

\newcommand{\sone}{S^1}

\newcommand{\stwo}{S^2}

\newcommand{\rone}{{\bf R}}
\newcommand{\zet}{{\bf Z}}

\newcommand{\del}{\partial}
\newcommand{\krondel}{\delta^{\mtiny{\text{Kronecker}}}}
\newcommand{\Hbar}{\overset{\,\_\!\_\!@!@!@!\_}{H}}
\newcommand{\pairing}{\langle@,@,@,,@!@!@!@!@!@!\rangle}
\newcommand{\ori}{{\frak o}}
\newcommand{\sw}{n^\ori_X}

\newcommand{\sast}{\sigma_{\!\ast}}

\newcommand{\Orth}{{O}}
\newcommand{\Ospin}{\Orth^\prime}

\newcommand{\Ok}{\Orth_{@!@!k}}
\newcommand{\Ospk}{\Orth^\prime_k}

\newcommand{\adid}{\cdot\{\pm\id\}}
\newcommand{\adsig}{\cdot\{\sigma_{\!\ast},\id\}}

\newcommand{\Hom}{\operatorname{Hom}}
\newcommand{\im}{\operatorname{im}}
\newcommand{\tors}{\operatorname{Tors}}
\newcommand{\rk}{\operatorname{rk}}
\newcommand{\id}{\operatorname{id}}
\newcommand{\Diff}{\operatorname{Diff}}

\newcommand{\Spinc}{\operatorname{Spin}_4^{c}}

\renewcommand{\a}{\alpha}
\renewcommand{\b}{\beta}
\renewcommand{\c}{\gamma}
\renewcommand{\d}{\delta}
\newcommand{\e}{\varepsilon}
\renewcommand{\l}{\lambda}
\newcommand{\s}{\sigma}
\newcommand{\z}{\zeta}

\newcommand{\tto}{\longrightarrow}
\newcommand{\inj}{\hookrightarrow}
\newcommand{\m}{\frac{\hspace{1.2mm}}{\hspace{1mm}}}

\newcommand{\X}{X_1}
\newcommand{\Xd}{X_d}
\newcommand{\Xq}{X_{1,q}}
\newcommand{\Xdq}{X_{d,q}}
\newcommand{\mult}{m_1,\ldots,m_n}
\newcommand{\Xdm}{X_{d;\mult}}
\newcommand{\Xdqm}{X_{d,q;\mult}}
\newcommand{\Xast}{X_\ast}

\newcommand{\Sq}{\Sigma_q}
\newcommand{\Sqq}{\Sigma_{q-1}}
\newcommand{\Sm}{S_{\mult}}
\newcommand{\Sast}{S_\ast}
\newcommand{\Sei}{S^3_m}

\begin{document}

\title{On the diffeomorphism groups of elliptic surfaces}
\author{Michael L\"onne}
\date{}
\maketitle

\begin{abstract}
\noindent
In this paper we determine for relatively minimal elliptic
surfaces with positive Euler number the image of the natural representation
of the group of orientation preserving  self-diffeomorphisms on $\Hbar$,
the second homology group reduced modulo torsion. To this end we construct
as many embedded spheres of square $\m 2$ such that an isometry not induced
from any combination of reflections at such spheres or from 'complex
conjugation' can be shown not to be induced from some diffeomorphism at
all. This is done with the help of Seiberg-Witten-invariants.
\end{abstract}

\vfill\vfill
\noindent
{\sc Introduction}\\[-.2cm]

\noindent
It is known that in dimension four corresponding notions of topological
and smooth manifold theory differ much more than in any other dimension.
This fact is best illustrated by the complete classification of simply
connected topological manifolds as given by Freedman \cite{f} versus
non-existence results for smooth structures on some topological manifolds
and multitudes of them on others which had been obtained by Donaldson
\cite{d} and others, e.g. in \cite{fm1}.\\ A further example is provided by
the respective groups of self-equivalences, homeo\-morphisms and
diffeomorphisms, which we assume to be orientation preserving without
further mentioning. They are studied via their representations on the
second homology and general results have been obtained on simply connected
manifolds:\\ The representations descend to faithful representations of
isotopy classes of homeomorphisms \cite{q},and pseudo-isotopy classes of
diffeomorphisms respectively \cite{k}. The group of homeomorphisms in this
case maps onto the orthogonal group \cite{q}, but the image of the
diffeomorphism group defies computation except for special manifolds
\cite{fm1},\cite{eo}, and after stabilization \cite{wall}.\\
It is the purpose of this paper to enlarge the set of special manifolds by
proving the following result, of which a special case answers a question
raised by Friedman and Morgan in their recent preprint \cite{fm2}:

\vfill
\pagebreak

\begin{main}
Let $X$ be a minimal elliptic surface with positive Euler number.
$\Diff(X)$ its group of diffeomorphisms and $\Hbar_2$ its second homology
reduced modulo torsion.
Then there exists an isometry $\sast$ such that
the image of the natural representation of $\,\Diff$ in the
orthogonal group $\Orth$ of $\Hbar_2$ with respect to the intersection
form is\\[.2cm]
\begin{tabular}{lp{10.5cm}}
$\quad\Ok\adid$ & the subgroup generated by $\pm\id$ and the
stabilizer group of the homology class of the canonical divisor in
the case of non-rational surfaces with geometric genus zero,
\\[.1cm]
$\quad\Ospk\adsig$ & the subgroup generated by $\sast$ and
the stabilizer group of the canonical class in the
subgroup of elements with real spinor norm one in the case of surfaces with
positive geometric genus,\\[.1cm]
$\quad\Orth$ & in the case of rational surfaces.
\end{tabular}
\end{main}


\noindent
We will address the existence and non-existence claims of the theorem
separately.
So in the first part we prove the existence of 'many' smoothly embedded
$2$-spheres of (intersection-)square $\m 2$ and exploit the interplay
between these spheres, the diffeomorphisms they give rise to, and algebraic
results on groups of orthogonal transformations.
This is just the sort of argument Friedman and Morgan set up when proving
the existence result for regular elliptic surfaces $\Xd$ without multiple
fibres and of positive Euler number $e(\Xd)=12d>0$ \cite{fm}. They also got
corollaries for regular elliptic surfaces $\Xdm$ with fibres of
multiplicities $\mult$.\\
Based on their first result we will be able to generalize it to arbitrary
elliptic surfaces $\Xdq$ without multiple fibres, of positive Euler number
$12d$ and irregularity $q$ possibly non-vanishing. These surfaces fibre
over curves which are smooth closed $2$-manifolds $\Sq$ of genus $q$. We
will then proceed to the case of elliptic surfaces $\Xdqm$ obtained from
the former ones by logarithmic transformations, i.e. by inserting multiple
fibres. In the regular case we improve the result of Friedman and Morgan
by a small margin.\\
Our result just fails to cover all minimal elliptic surfaces since we
restrict ourselves to the case of positive Euler number. The remaining
surfaces will require essentially different methods because even in the
topological category only few homology classes are represented by
spheres.\\
The main tool in the second part are Seiberg Witten invariants which have
already proved invaluable in the realms of smooth four manifolds. We start
setting up the terminology in line with the definitions in \cite{b}. This
provides the means to link the results from \cite{b} to the applications
we have in mind. Our non-existence claims are then proved using the
invariance of the homology class of a canonical divisor and exploiting the
subtleties of the homology orientation very much in the spirit of
\cite{eo}.\\
Taking into account the existence of a diffeomorphism induced from complex
conjugation, cf. \cite{fm1}, is the final bit to accomplish the proof of
the theorem.\\
In the course of the argument we will freely use Poincar\'{e} duality to
identify homology classes of divisors with the Chern classes of the line
bundles they determine. A fair knowledge of the geometry of elliptic
surfaces as can be found in \cite{fm} or \cite{bpv} is assumed
tacitly.

\pagebreak

\noindent
{\sc Embedded spheres in elliptic surfaces}\\[-.2cm]

\begin{thm}[\cite{fm}]
\label{X_spheres}
Every class in $H_2(\X-nf)$ of square $\m 2$ is represented by a $2$-sphere
smoothly embedded in $\X-nf$. In fact there is a set of $9+2n$ such spheres
$$
\{\a_1,\ldots,\a_{2n+1},\b_1,\ldots,\b_8\}
$$
whose classes form a basis for $H_2(\X-nf)$ such that
\begin{itemize}
\item the $\{\b_j\}$ intersect according to the Dynkin diagram of $\m E_8$,
\item the intersection pattern of the $\{\a_i\}$ is given by a complete
graph on $2n+1$ vertices whose edges are weighted by the intersection
number $\m 2$,
\item every $\a_i$ intersects $\b_8$ algebraically once,
\item $\{\ell+\a_i\}$ represents any given basis for the radical of the
intersection pairing on $H_2(\X-nf)$; $\ell$ denoting the long vector of
$\m E_8$.
\end{itemize}
Moreover a set $\{\b_j^\prime\}$ with the first property can be
supplemented to such a basis.
\end{thm}

\proof
Except for the addendum the claim is that of thm.5.1 in
\cite{fm} generalized along the lines of thm.5.11, loc.cit..\\
Given any basis as in the theorem we get another for every isometry of
$(H_2,q_H)$. Such an isometry we get by mapping $\{\b_j\}$ of the theorem to
the $\{\b_j^\prime\}$ and extending by an isomorphism of the radical. The
images of the $\{\a_i\}$ then are represented by spheres according to the
first assertion of the theorem and are easily seen to supplement the
$\{\b_j^\prime\}$ to a basis.
\qed

\begin{lemma}
\label{Xdq_decomposition}
An elliptic surface $\Xdq$ which is properly irregular, i.e. $q\geq1$, and
has positive Euler number $12d$ decomposes as
$$
\Xdq=\underbrace{\X \#_f\ldots \#_f\X}_{d}\#_{2f}B
$$
with $B:=\Sqq\times f$. This is equivalent to
$$
\Xdq-f=\underbrace{(\X-2f)\cup_{f\times\sone}\ldots
\cup_{f\times\sone}\!(\X-2f)}_{d-1}
\cup_{f\times\sone}(\X-3f)\cup_{2(f\times\sone)}\!(B-2f).
$$
\end{lemma}

\proof
There is a well known decomposition, cf. \cite[pp. 159, 190, 195]{fm}:
$$
\Xdq=\underbrace{\X\#_f\ldots \#_f\X}_{d}\#_fB^\prime,
$$
with $B^\prime=\Sq\times f$ which can be considered as lifted from a
decomposition of the base
$$
\Sq=\underbrace{\stwo \#_{pt}\ldots \#_{pt}\stwo}_{d}\#_{pt}\Sq.
$$
Since $\Sq$ is $\Sqq\#_{2pts}\stwo$ and we may transfer the new $\stwo$
summand to the neighbouring one, we get:
$$
\Sq=\underbrace{\stwo \#_{pt}\ldots \#_{pt}\stwo}_{d-1}\#_{pt}
(\stwo\#_{pt}\stwo)\#_{2pts}\Sqq,
$$
which induces the claimed decomposition on the total space.
Cut out a regular fibre and reformulate the fibre connected sum in terms of
the union along submanifolds to get the second assertion.
\qed

\begin{lemma}
\label{cylinder}
There is a basis of $H_2(B,2f)$ geometrically represented by
\begin{itemize}
\item a section to the fibration map,
\item smoothly embedded cylinders with vanishing selfintersection and the
distinct boundary components mapping to the distinct fibres.
\end{itemize}
\end{lemma}

\proof
$(B,2f)=(f\times(\Sqq,2pts))$. Choose injective maps
\mbox{$\psi_i:\sone\inj f$} and \mbox{$\phi_j:(I,\{0,1\})\inj(\Sqq,2pts)$}
representing bases of $H_1(f),\,H_1(\Sqq,2pts)$ respectively. The product
maps
$$
\Psi_{i,j}:(\sone\times I,\sone\times\{0,1\})\inj(B,2f)
$$
then represent classes in $H_2(B,2f)$. These are obviously embedded
cylinders of the kind we want and a look at the exact homology sequence of
the pair $(B,2f)$ should convince us that together with a section
they form a basis for $H_2(B,2f)$
\qed

\begin{lemma}
\label{thimble}
There are four smoothly embedded vanishing cells of square $\m 1$ disjoint
off their boundaries
$$
(D^2,\sone)\inj(\X-f,2f)
$$
representing relative classes in $H_2(\X-f,2f)$ such that their boundaries
are arbitrarily given simple closed curves generating $H_1(2f)$.
\end{lemma}

\proof
Without loss of generality we may assume that $\X-f$ contains six cusp
neighbourhoods disjoint from $2f$ and each other. Given a curve then we
construct a cell as follows:
\pagebreak\\
\noindent
Choose a regular fibre in the boundary of the cusp neighbourhood. This
fibre and the given curve map to points in the base which we connect by a
simple arc such that the fibration over it is trivial. So we may transfer
the homology class of the curve to the other end. Pick an associated
vanishing cell. Its boundary can be isotoped over the arc to the given
curve. Glue the cylinder thus obtained to the cell to get the
desired embedding. Different cells won't intersect as long as we choose
different cusp neighbourhoods for different curves and disjoint arcs -
which obviously is a condition easily met.
\qed

\begin{thm}
\label{Xq_spheres}
In the proper irregular case, i.e. $q\geq1$,
there is a generating set for $H_2(\Xq-f)$ of square $\m 2$ classes which
are represented by $2$-spheres
$$
\{\a_1,\ldots,\a_7,\b_1,\ldots,\b_8,\c_1,\ldots,\c_{2q-1},
\d_1,\ldots,\d_{2q-1},\e_1,\e_2,\e_3\},
$$
smoothly embedded in $\Xq-f$ and such that
\begin{itemize}
\item the $\{\a_i\}$ intersect according to a complete graph on 7
vertices: $\a_i\cdot\a_{i^\prime}=-2$,
\item the $\{\b_j\}$ intersect according to the Dynkin diagram of $\m E_8$,
\item every $\a_i$ intersects $\b_8$ algebraically once,
\item the intersection graph spanned by
$\{\b_8,\c_1,\d_1,\e_1,\e_2,\e_3\}$ is:\\
\\
\unitlength2.5cm
\begin{picture}(4.5,2.2)(.1,.9)
\put(1,2){\circle*{.06}}
\put(1.03,1.9){\mathbbox{\c_1}}
\put(1,2){\line(1,0){3}}
\put(1.5,2.1){\mathbox{1}}
\put(2.7,2.1){\mathbox{-2}}
\put(3.5,2.1){\mathbox{1}}
\put(2,2){\circle*{.06}}
\put(1.96,1.9){\mathbbox{\e_1}}
\put(2,2){\line(1,2){.5}}
\put(2.1,2.5){\mathbox{-2}}
\put(2,2){\line(1,-2){.5}}
\put(2.1,1.5){\mathbox{1}}
\put(3,2){\circle*{.06}}
\put(3.06,1.9){\mathbbox{\e_2}}
\put(3,2){\line(-1,2){.5}}
\put(2.9,2.5){\mathbox{-2}}
\put(3,2){\line(-1,-2){.5}}
\put(2.9,1.5){\mathbox{1}}
\put(4,2){\circle*{.06}}
\put(4.03,1.9){\mathbbox{\d_1}}
\put(2.5,3){\circle*{.06}}
\put(2.53,3.1){\mathbbox{\e_3}}
\put(2.5,3){\line(0,-1){.95}}
\put(2.5,1){\circle*{.06}}
\put(2.53,.89){\mathbbox{\b_8}}
\put(2.5,1){\line(0,1){.95}}
\put(2.6,1.6){\mathbox{1}}

\end{picture}
\end{itemize}
\end{thm}

\proof
Decompose $\Xq-f$ according to lemma 1:
$$
\Xq-f=(\X-3f)\cup_{2(f\times\sone)}(B-2f)=(\X-f)\#_{2f}B.
$$
This is a fibre connected sum along two disjoint regular fibres. Choosing
injective maps $\phi_{1,2}:\sone\inj f$ generating the first homology of
one of these we get a second pair of such maps to the other fibre by the
product structure on $B$. Together they give rise to embedded cylinders in
$(B,2f)$ and embedded cells in $(\X-f,2f)$ according to the preceding
lemmas.
\pagebreak
\\
\noindent
If we glue these in the obvious way we get smoothly embedded
$2$-spheres in $\Xq-f$ with square $\m 2$. Call them $\c_k,\d_l$ such that
the index refers to an indexing of a basis for $H_1(\Sqq,2pts.)$, cf. lemma
2, and the $\c$'s are those constructed from $\phi_1$, the $\d$'s those
from $\phi_2$.\\
Let's now have a closer look at the homology sequence of the pair
$(\Xq-f,\X-3f)$:
$$
H_2(\X-3f)\stackrel{i}{\to}H_2(\Xq-f)\stackrel{j}{\to}H_2(\Xq-f,\X-3f)
\stackrel{\del}{\to}H_1(\X-3f).
$$
Here replace the relative group by $H_2(B,2f)$ by means of a suitable
excision lemma, cf. \cite[\S3.0]{la}. Since the section in
$H_2(B,2f)$ maps neither to zero nor to a torsion class in $H_1(\X-3f)$
the image of the relative map is spanned by the cylinders of lemma
\ref{cylinder}. Of course our new spheres $\c_k,\d_l$ are lifts of those.
Thus any generating set for $H_2(\X-3f)$ will be supplemented by these
spheres to a generating set for $H_2(\Xq-f)$.\\
Pick such a set $\{\a_i^\prime,\b_j^\prime\}$ as given by theorem
\ref{X_spheres} such that the basis of the radical
$\{\ell^\prime+\a_i^\prime\}$ is dual to $\c_k,\d_l$, i.e.:
$$
\langle\ell^\prime+\a_i^\prime,\c_k\rangle=\krondel_{i,1},\quad
\langle\ell^\prime+\a_i^\prime,\d_l\rangle=\krondel_{i,2}
$$
Note that we may speak unambiguously of intersection with $\c_k$, for it is
determined by the relative class it represents in $H_2(\X-f,2f)$ which is
independent of the index, ditto for $\d_l$.\\
Now we adjust this set for our purposes. First set:
$$
\b_i:=\b_i^\prime\quad i=1,\ldots,7;\qquad\b_8:=\b_8^\prime-\ell^\prime-
\langle\b_8^\prime-\ell^\prime,\c_k\rangle\a_1^\prime-
\langle\b_8^\prime-\ell^\prime,\d_l\rangle\a_2^\prime.
$$
It is easily checked that:
\begin{itemize}
\item the $\{\b_j\}$ span an intersection graph $\m E_8$,
\item $\b_8\cdot\c_k=0=\b_8\cdot\d_l$.
\end{itemize}
With theorem \ref{X_spheres} we may supplement the new $\b_j$ by a set
$\{\a_i\}$  subject to a duality property as above relative to the new
$\ell$.\\
Finally define the $\e_1,\e_2,\e_3$ to be :
\begin{eqnarray*}
& \e_1=\langle\c_k,\ell\rangle(\ell+\a_1)+
 \langle\d_l,\ell\rangle(\ell+\a_2)+\a_1 \\
& \e_2=\langle\c_k,\ell\rangle(\ell+\a_1)+
 \langle\d_l,\ell\rangle(\ell+\a_2)+\a_2 \\
&\e_3=\langle\c_k,\ell\rangle(\ell+\a_1)+
 \langle\d_l,\ell\rangle(\ell+\a_2)+\a_3
\end{eqnarray*}
And now we are finished. The intersection pattern of the set chosen is
seen to meet the requirements of the theorem. Moreover the set has already
been constructed as a set of spheres except for $\e_1,\e_2,\e_3$, but
as they can be considered as elements in $H_2(\X-3f)$ there are embedded
spheres representing them by theorem \ref{X_spheres}.
\qed

\begin{cor}
\label{Xdq_spheres}
In the proper irregular case there is a generating set for $H_2(\Xdq-f)$ of
square $\m 2$ classes which are represented by $2$-spheres
$$
\left\{\a^1_{1,...,7},\b_j^1,\c_k,\d_l,\e_{1,2,3},
\a_i^n,\b_j^n,\z_{1,2}^n|
{\mtiny\renewcommand{\arraystretch}{.6}
\begin{array}{r@{}c@{}l@{\mtiny,\,}r@{}c@{}l}
 	 \mtiny 1\leq & \mtiny i & \mtiny \leq 5 &
     \mtiny 1\leq & \mtiny k,l & \mtiny \leq 2q-1\\
   \mtiny 1\leq & \mtiny j & \mtiny \leq 8 &
     \mtiny 2\leq & \mtiny n & \mtiny \leq d\\
\end{array}}
\right\}
$$
smoothly embedded in $\Xdq-f$ and such that:
\begin{itemize}
\item the $\{\a_i^1,\b_j^1,\c_k,\d_l,\e_{1,2,3}\}$ intersect as in theorem
      \ref{Xq_spheres},
\item the $\{\a_i^n,\b_j^n\}$ intersect as in theorem \ref{X_spheres},
\item for every sphere $\z_{1,2}^n$ there are two spheres among the $\a$'s
      of the $n^{th}$ and the $(n$-1$)^{th}$ summand respectively which
      intersect $\z$ algebraically once.
\end{itemize}
\end{cor}

\proof
Start with a decomposition of $\Xdq-f$ as given by lemma
\ref{Xdq_decomposition}:
$$
\Xdq-f=\Xq\#_f\X\ldots\#_f\X\#_f(\X-f).
$$
By a Mayer-Vietoris argument a generating set is given as soon as there
are such sets for each summand and spheres which after restriction
generate the first homology of all those fibres on which fibre connected
sum has been performed.\\
Sets of the first kind we get from theorem \ref{X_spheres} and theorem
\ref{Xq_spheres} so our set will automatically be furnished with the first
two properties.\\
For the construction of the spheres of the second kind recall that two
fibres and their neighbourhoods in the summands had to be identified in
order to perform fibre connected sum. Choose generating curves for the first
homology of these fibres compatible with the identification. Then inside
the connected sum we can glue the vanishing cells in each summand which we
obtain as in lemma \ref{thimble}. The outcome is obviously a $2$-sphere
embedded in the connected sum of the adjacent summands of the fibre we
started with.\\ By construction in each of this summands a sphere $\z$  can
be considered as a relative homology class, thus there is a class $\ell+\a$
dual to it, cf. the proof of the preceding theorem. Hence
$$
\langle\z,\langle\z,\ell\rangle(\ell+\a)+\a\rangle=1,
$$
and we are left with showing that we may even assume the class
$\langle\z,\ell\rangle(\ell+\a)+\a$ to be represented by one of the
$\a_i$'s. But this is indeed the case as its square is $\m 2$ and the sum
with $\ell$ is an element of the radical for the summand under
consideration, thus a isometry as in the proof of theorem \ref{X_spheres}
will do the job.
\qed

\begin{thm}
\label{Xdq_diffgrp}
If $\Xdq$ is an elliptic surface without multiple fibres and Euler number
$12d>0$ and $f$ is a generic regular fibre, then:
\begin{itemize}
\item every class in $H_2(\Xdq-f)$ of square $\m 2$ is represented by a
$2$-sphere smoothly embedded in $\Xdq-f$,
\item in the rational case, i.e. $\Xdq\!=\!\X$, reflections in
classes of square $\m 2$ generate the group of all isometries of
$H_2(\Xdq-f)$ which restrict to the identity on the radical, in all other
cases they generate a subgroup of index at most $2$,
\item every such reflection is realized by a diffeomorphism which is the
identity on a neighbourhood of $f$.
\end{itemize}
\end{thm}

The proof mimics that given for thm.6.2 in \cite{fm}. Its geometric
ingredient is the fact that a $2$-sphere of square $\m 2$ gives rise to a
diffeomorphism which is the identity outside a neighbourhood of the sphere
and which acts on homology as reflection in the class of the sphere, cf.
\cite[p.358f]{fm1}. The algebraic counterpart is a theorem due to Ebeling
in the general case and the slightly enhanced lemma 5.9 of \cite{fm} in the
rational case:

\begin{thm}[\cite{e}]
\label{ebeling}
Let $(L,\langle,\rangle)$ be an even lattice. Let $\Delta\subset L$ be a
set of vectors of square $\m 2$, and let $\Gamma_\Delta$ be the subgroup of
the isometry group of $L$ generated by reflections in $\l\in\Delta$.
Suppose that $\Delta$ satisfies the conditions i), ii), and iii) below:
\begin{enumerate}
\item $\Delta$ spans $L$,
\item $\Delta$ is contained in a single $\Gamma_\Delta$-orbit,
\item in $\Delta$ there are six elements $\l_1, ...,\l_6$ which intersect
as given by the following diagram:\\
\unitlength2.5cm
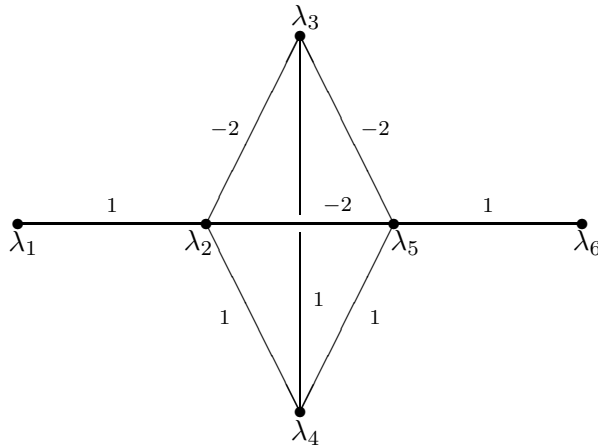
\begin{picture}(4.7,2.2)(.1,.9)
\put(1,2){\circle*{.06}}
\put(1.03,1.9){\mathbbox{\l_1}}
\put(1,2){\line(1,0){3}}
\put(1.5,2.1){\mathbox{1}}
\put(2.7,2.1){\mathbox{-2}}
\put(3.5,2.1){\mathbox{1}}
\put(2,2){\circle*{.06}}
\put(1.96,1.9){\mathbbox{\l_2}}
\put(2,2){\line(1,2){.5}}
\put(2.1,2.5){\mathbox{-2}}
\put(2,2){\line(1,-2){.5}}
\put(2.1,1.5){\mathbox{1}}
\put(3,2){\circle*{.06}}
\put(3.06,1.9){\mathbbox{\l_5}}
\put(3,2){\line(-1,2){.5}}
\put(2.9,2.5){\mathbox{-2}}
\put(3,2){\line(-1,-2){.5}}
\put(2.9,1.5){\mathbox{1}}
\put(4,2){\circle*{.06}}
\put(4.03,1.9){\mathbbox{\l_6}}
\put(2.5,3){\circle*{.06}}
\put(2.53,3.1){\mathbbox{\l_3}}
\put(2.5,3){\line(0,-1){.95}}
\put(2.5,1){\circle*{.06}}
\put(2.53,.89){\mathbbox{\l_4}}
\put(2.5,1){\line(0,1){.95}}
\put(2.6,1.6){\mathbox{1}}

\end{picture}
\end{enumerate} Then $\Gamma_\Delta$ is the subgroup of the isometry group
of $(L,\langle,\rangle)$ consisting of automorphisms of real spinor norm
one which are the identity on $\Hom(L,\zet)/\im L$, and
$\Gamma_\Delta\cdot\Delta$ is the set of all vectors of square $\m 2$ in $L$.
\qed
\end{thm}

\begin{lemma}
\label{semi_definite}
Let $(L,\langle,\rangle)$ be a lattice which decomposes as a direct sum of
the radical and a unimodular lattice $E$. Let $\Delta\subset L$ be a set
of vectors of square $\m 2$. Suppose that:
\begin{enumerate}
\item $\Delta$ spans $L$,
\item $\Gamma_{\Delta\cap E}\cdot(\Delta\cap E)$ is a single
$\Gamma_{\Delta\cap E}$-orbit and is the set of all vectors of square $\m 2$
in $E$,
\item $\Gamma_{\Delta\cap E}$ contains all isometries of
$(E,\langle,\rangle_E)$.
\end{enumerate}
Then $\Delta$ is contained in a single $\Gamma_\Delta$-orbit which
consists of all vectors of square $\m 2$ in $L$ and $\Gamma_\Delta$ is the
group of all isometries of $(L,\langle,\rangle)$ which are the identity on
the radical.
\end{lemma}


{\it Proof of theorem \ref{Xdq_diffgrp}: }
As the theorem has been proved in \cite{fm} in the regular case we may
first consider irregular surfaces; so we may take $\Delta$ to be the
set of spheres given in the corollary to theorem \ref{Xq_spheres}. Thereby
the hypotheses $i)$ and $iii)$ of theorem \ref{ebeling} are automatically
satisfied. To check hypothesis $ii)$ we make use of the following
criterion: two elements of $\Delta$ are in the same orbit of
$\Gamma_\Delta$ if they intersect algebraically once up to sign. Since
being conjugate is a transitive relation the verification of hypothesis
$ii)$ boils down to checking whether the intersection graph is connected by
edges which carry intersection numbers $\pm1$. This in fact is ensured by
the various intersection properties of the preceding results.\\
In the rational case we take $\Delta$ to be the set of spheres given in
thm.\ref{X_spheres}. Then the first hypothesis of the lemma is satisfied,
and so are the others by well known properties of the $E_8$ lattice.\\
Now let's exploit the conclusions: In both cases $\Gamma_\Delta$ is shown
to be the group of all reflections in classes of square $\m 2$. Its
generators are given by spheres which are disjoint from the fibre $f$. By
the geometrical ingredient mentioned above the generators are realized by
diffeomorphisms which are the identity on a neighbourhood of
$f$, and so is the entire group! Moreover because there is only one orbit
for classes of square $\m 2$, an orbit of an embedded sphere under the
diffeomorphisms just constructed, will supply embedded spheres representing
all classes of square $\m 2$.\\
Finally note that in the general case the homology lattice modulo its
radical is unimodular. Hence we may identify the group of isometries which
induce the identity on $\Hom(L,\zet)/\im L$ with the group of isometries
which are the identity on the radical. Thus of the latter $\Gamma_\Delta$
is a subgroup of index at most two, since the real spinor norm is a group
homomorphism to $\zet_2$. In the rational case the stronger conclusion is
already given by the lemma.
{\hfill $\Box$ \vspace{.6cm}}

\begin{cor}
If $\Xdq$ is an elliptic surface without multiple fibres and Euler number
$12d>0$ and $2f$ are two disjoint regular fibres, then:
\begin{itemize}
\item every class in $H_2(\Xdq-2f)$ of square $\m 2$ is represented by a
$2$-sphere smoothly embedded in $\Xdq-2f$,
\item two such spheres are conjugate under the group of reflections on
classes of square $\m 2$.
\end{itemize}
\end{cor}
\hfill $\Box$

\pagebreak

\begin{lemma}
\label{Xdqm_decomposition}
An elliptic surface $\Xdqm$ with positive Euler number has a decomposition
$$
\Xdqm-f=(\Xdq-2f)\cup_{f\times\sone}(\Sm-f)=(\Xdq-f)\#_f\Sm,
$$
with $\Sm$ a Seifert fibration of tori over $\stwo$ such that its multiple
fibres have multiplicities $m_i$.\\
This decomposition induces an exact homology sequence:
$$
H_i(\Xdq-2f)\to H_i(\Xdqm-f)\to H_i(\Sm,f)\to H_{i-1}(\Xdq-2f).
$$
\end{lemma}

\proof
A decomposition of the base curve of the fibration into a disc to which
all multiple fibres map and the remaining part to which all the other
singular fibres map lifts to the decomposition of the total space as in the
claim above.\\
The exact sequence is then derived from the exact homology sequence of the
pair $(\Xdqm,\Xdq-2f)$ and the isomorphism, cf. \cite[\S3.0]{la}:
$$
\Psi:H_i(\Xdqm-f,\Xdq-2f)\stackrel{\simeq}{\to}H_i(\Sm,f).
$$
\qed

\begin{lemma}
\label{fibre_detect}
There is an exact sequence
$$
H_2(\Xdq-2f)\to\Hbar_2(\Xdqm-f)\to\zet_m\to0,
$$
with $\Hbar$ the quotient group of $H$ modulo its subgroup of torsion
elements.\\
In this sequence the multiple fibres taken as geometric preimages under the
fibration map, i.e. without their multiplicities, represent classes which
map to a set of generators for the quotient group $\zet_m$ and $m$ is
given as the least common multiple of the multiplicities.
\end{lemma}

\noindent
To make the following argument clearer at least in a typographical respect
we introduce $\tors_i$ as a shorthand notation for the torsion subgroup of
$H_i$, as well as:
$$
X:=\Xdq,\quad\Sast:=\Sm,\quad\Xast:=\Xdqm.
$$

{\it Proof of lemma \ref{fibre_detect}:}
Take the following part of the sequence from the lemma above and compute
the rank for each group:
$$
\underset{\rk=1}{H_2(\Sast,f)}@>{\del\circ\Psi^{-1}}>>
\underset{\rk=2q+1}{H_1(X-2f)}\to\underset{\rk=2q}{H_1(\Xast-f)}
\to\underset{\rk=0}{H_1(\Sast,f)}.
$$
Thus the rank of $\ker\del\circ\Psi^{-1}$ vanishes. Moreover $H_1(X-2f)$
is torsionfree, so in fact $\ker\del\circ\Psi^{-1}=\tors_2(\Sast,f)$ and
we can write down an exact sequence:
$$
H_2(X-2f)\to H_2(\Xast-f)\to\tors_2(\Sast,f)\to0.
$$
Since $H_2(X-2f)$ is torsionfree, the torsion subgroup of $H_2(\Xast-f)$
is mapped isomorphically onto its image. Dividing both out yields:
$$
H_2(X-2f)\stackrel{j}{\to}\Hbar_2(\Xast-f)\to
\tors_2(\Sast,f)_{\displaystyle /\tors_2(\Xast-f)}\to0.
$$
This must be shown to be the sequence of the claim. To this end consider
the following commutative diagram of exact rows and columns:\\
$$
\begin{array}{ccccl}
 & & & & H_2(f)\\
 & & & \swarrow & \quad\downarrow\\
H_3(\stwo\times f,f) & \tto & H_2(X-2f) & \tto & H_2(X-f)\\
\downarrow & & \downarrow {\scriptstyle j} & & \quad\downarrow\\
H_2(\Sast-f) & \tto & H_2(\Xast-f) & \tto & H_2(X-f,f)\\
\end{array}
$$
The left hand column is part of the homology sequence of the pair
$(\Sast,\Sast-f)$ modified by the excision lemma from \cite{la}.
With duality and computations of fundamental groups we obtain ranks for
the following part of this sequence:
$$
\underset{\rk=0}{H_3(\Sast-f)}\to\underset{\rk=2}{H_3(\Sast)}
\to\underset{\rk=2}{H_3(\stwo\times f,f)}\to H_2(\Sast-f).
$$
thus the last map which is the map of the diagram maps to
$\tors_2(\Sast-f)$. Keeping this in mind a diagram chase shows that a class
in the image of $j$ and commensurable with $[f]$ in $H_2(\Xast-f)$ must be
an integer multiple of $[f]$.\\
On the other hand a multiple fibre without its multiplicity represents the
class $\frac{1}{m_i}[f]$ hence $[f]$ is at least $m$-divisible in
$H_2(\Xast-f)$ with $m=lcm(m_i)$. To state it differently the
multiple fibres without their multiplicities represent classes which
generate a subgroup $\zet_m$ in $\tors_2(\Sast,f)/\tors_2(\Xast-f)$.\\
Computation using duality and fundamental groups
$$
\begin{array}{lllll}
\tors_2(\Sast,f) & \cong & \tors_1(\Sast-f) & \cong & \bigoplus
\zet_{m_i}\\
\tors_2(\Xast-f) & \cong & \tors_1(\Xast,f) & \cong & \bigoplus
{\zet_{m_i}}_{\displaystyle/\zet_m}\\
\end{array}
$$
proves this inclusion to be an isomorphism since the cardinalities are
shown to coincide, and so we are done.
\qed

\begin{thm}
\label{Xdqm_spheres}
Let $mf_m$ be a multiple fibre of multiplicity $m$ in the elliptic surface
$\Xdqm$, so $m=m_i$ for some $1\leq i\leq n$.\\
Then in the complement of a generic regular fibre there exist three
smoothly embedded spheres $\s_{1,2,3}$ of square $\m 2$, such that:
\begin{itemize}
\item $\s_3$ is disjoint from all multiple fibres, $\s_{1,2}$ from all
multiple fibres except $mf_m$,
\item in $H_2(\Xdqm)$ the difference of the first two spheres represents
$[f_m]$,
\item $\s_1\cdot\s_3=\s_2\cdot\s_3=\pm1$,
\item the class of $\s_3$ is in the image of the map
$$
H_2(\Xdq-2f)\to H_2(\Xdqm-f).
$$
\end{itemize}
\end{thm}

\proof
A fibred neighbourhood of a multiple fibre is diffeomorphic to
$\sone\times\Sei$, the product of the loop with a Seifert fibration of
loops over the disk with a central fibre of multiplicity $m>1$. As
sometimes earlier we may assume that there exist three cusp
neighbourhoods.\\
A path
$$
w:(I,\del I)\tto(\Sei,\del\Sei)
$$
which meets the central fibre in exactly one point gives rise to a
cylinder of square $\m 2$ in the fibred neighbourhood simply by taking the
cartesian product with the loop. By methods exploited already in lemma
\ref{thimble} we may glue two vanishing cells from different cusp
neighbourhoods to this cylinder to get the embedded sphere $\s_2$.\\
Apply the same construction to a new path. Up to the central fibre it
coincides with the old one, then pursues the fibre once before taking up
the old path again up to its end. Just homotop it off itself and the
construction yields another embedded sphere $\s_1$.\\
The sphere $\s_3$ is obtained by gluing one of the vanishing cells just
used to another one in the third cusp neighbourhood. These spheres can now
be seen to satisfy the claims of the theorem:
\begin{itemize}
\item the neighbourhoods we used can be chosen disjoint from all the other
multiple fibres and a generic regular one; out of these neighbourhoods the
constructed spheres map to simple arcs on the base and these again can be
chosen in an appropriate way,
\item the difference of the classes represented by the first two spheres
is represented by the cartesian product of the central fibre in the
Seifert fibration with the loop,
\item the intersection with the third sphere is by construction up
to sign just the selfintersection of the vanishing cell in the first cusp
neighbourhood,
\item because of the first assertion the third sphere is liftable.
\end{itemize}
\qed

\begin{thm}
\label{Xdqm_diffgrp}
If $\Xdqm$ is an elliptic surface with non-vanishing Euler number and $f$
a generic regular fibre, then:
\begin{itemize}
\item every class in $\Hbar_2(\Xdqm-f)$ of square $\m 2$ is represented by a
$2$-sphere smoothly embedded in $\Xdqm-f$,
\item in the case of surfaces with $p_g=0$, reflections in
classes of square $\m 2$ generate the group of all isometries of
$\Hbar_2(\Xdqm-f)$ which restrict to the identity on the radical, in all
other cases they generate the subgroup of index at most $2$ which consists
of all element of real spinor norm one,
\item every such reflection is realized by a diffeomorphism which is the
identity on a neighbourhood of $f$.
\end{itemize}
\end{thm}

\proof
The claim is proved already for elliptic surfaces without multiple fibres.
Moreover we are left to check the hypotheses for Ebeling's theorem and
lemma \ref{semi_definite} in the appropriate cases since the rest of the
proof of theorem \ref{Xdq_diffgrp} goes through unchanged. But they are a
consequence of the exact sequence of lemma
\ref{fibre_detect}
$$
H_2(\Xdq-2f)\to\Hbar_2(\Xdqm-f)\to\zet_m\to0,
$$
and the assertion proved there that the span of classes represented by the
multiple fibres maps surjective onto the quotient:
\begin{itemize}
\item the image of $H_2(\Xdq-2f)$ has a generating set represented by
spheres according to the corollary to theorem \ref{Xq_spheres}; enlarge
this set by three spheres as in theorem \ref{Xdqm_spheres} for each
multiple fibre, then multiple fibres are given as differences of spheres,
thus the total lattice is generated,
\item the intersection numbers in theorem \ref{Xdqm_spheres} are such that
the new spheres are conjugate to spheres representing elements in
$H_2(\Xdq-2f)$, in this lattice all spheres of square $\m 2$ are conjugate
by the corollary to theorem \ref{Xdq_diffgrp},
\item since in $\Hbar_2(\Xdqm-f)$ there is a sublattice isomorphic to
$H_2(\Xdq-f)$, the special intersection lattice of the last hypothesis
in Ebeling's result can be found again in case that $p_g>0$,
\item by the same argument the lattice $\Hbar_2(\Xdqm-f)$ decomposes as a
direct sum of the radical and an $\m E_8$ lattice in the case that $p_g=0$.
\end{itemize}
\qed

\begin{cor}
Let $X$ be an elliptic surface with non-vanishing Euler number, then on
homology the diffeomorphisms just constructed generate the stabilizer group
of the class of a general fibre in case that $p_g(X)=0$, they generate the
elements of real spinor norm one only in case that $p_g(X)>0$.
\end{cor}

\proof
This is immediate from the fact that the isometries of the homology lattice
$H_2(X-f)$ stabilizing its radical map surjectively onto the isometries of
the lattice $H_2(X)$ stabilizing the class of a general fibre. This
argument is given in greater detail in the proof of thm.6.5 in \cite{fm}.
\qed

\begin{thm}
\label{E+K3}
Let $X$ be an Enriques surface or a $K3$ surface, then:
\begin{itemize}
\item
every class in $\;\Hbar_2(X)\,$ of square $\m 2$ is represented by a
$2$-sphere smoothly embedded in $X$,
\item
the reflections in classes of square $\m 2$ generate the group of all
isometries of $\Hbar_2(X)\,$ with real spinor norm one,
\item
every such reflection is realized by a diffeomorphism, so diffeomorphisms
of $X$ generate the group of isometries of real spinor norm one.
\end{itemize}
\end{thm}

\noindent
This extension to the previous results is possible since there are special
$K3$ surfaces and Enriques surfaces which by their complex structure have
sections and two-sections respectively which provide embedded spheres of
square $\m 2$ complementing the previous sets to a basis of $\Hbar_2$.
This fact has already been exploited in the $K3$ case in the proof of
thm.6.6 in \cite{fm}. Therefore and since the proofs are conceptionally
analogous we restrict ourselves to the case of Enriques surfaces. Again
we make use of an algebraic result:

\begin{lemma}[\cite{e}]
Let $L$ be a lattice and $\Delta\subset L$ be a set of vectors of square
$\m 2$ such that:
\begin{enumerate}
\item
$\Delta$ is a basis for $L$,
\item
the elements of $\Delta$ intersect according to the Dynkin diagram of
$\m E_{10}$.
\end{enumerate}
Then $\Gamma_\Delta$ is the subgroup of the isometry group of $L$
consisting of automorphisms of real spinor norm one and $\Delta$ is
contained in a single $\Gamma_\Delta$-orbit which consists of all vectors
of square $\m 2$ in $L$.
\end{lemma}


{\it Proof of theorem \ref{E+K3}:}
As already stated there is a sphere representing a class $\s$ of square
$\m 2$ such that $\Hbar_2(X)$ is spanned by $\s$ and classes from
$\Hbar_2(X-f)$. In the image of $\Hbar_2(X-f)$ choose a set of classes
$\{\b_i^\prime\}_{1\leq i\leq8}$ intersecting according to the Dynkin
diagram of $\m E_8$ and let $f$ be the unique class in the radical which
pairs to one with $\s$. Then set
$$
\begin{array}{lll}
\b_i & = & \b_i^\prime-\langle\b_i^\prime,\s\rangle f\qquad1\leq i\leq8\\
\b_9 & = & f-\ell\\
\b_{10} & = & \s,
\end{array}
$$
where $\ell$ is the long vector of the $\m E_8$ spanned by
$\{\b_i\}_{1\leq i\leq8}$, and it is to check these $\b_i$ intersect
according to the Dynkin diagram of $\m E_{10}$. Moreover they all are
represented by embedded spheres since the $\{\b_i\}_{1\leq i\leq9}$ are
given by classes from $\Hbar_2(X-f)$ which are all represented by spheres
according to thm.\ref{Xdqm_diffgrp}.\\
We finally apply the lemma and conclude as in the proof of
thm.\ref{Xdq_diffgrp}: Since generators of $\Gamma_\Delta$ are induced
from diffeomorphisms, so is the whole group of reflections which proves
part of the third claim. There are classes which are represented by
embedded spheres and hence they are all, for they form a single orbit
under $\Gamma_\Delta$, proving the first claim. The rest is immediate
from the lemma.
\hfill $\Box$

\pagebreak

\noindent
{\sc Applications of Seiberg-Witten invariants}\\[-.2cm]

\noindent
We introduce the following terminology for Seiberg-Witten
theory which is for convenience (almost) that of that of our main
reference, cf.
\cite[pp. 11,12,14]{b}.\\

\begin{tabular}{p{4cm}p{8.5cm}}
$SC$-structure & an equivalence class of a complex vector bundle $W$,
which is fibrewise irreducible as module of the Clifford bundle; the set
of $SC$-structures corresponds bijectively to the set of
$\Spinc$-structures.\\[.2cm]
\mbox{$SW$-multiplicity map} (case $b^+>1$) & the Seiberg-Witten map $n_X$
which maps a
$SC$-structure $W$ to the Euler class of a moduli space of solutions to
the monopole equations on $W$.\\[.2cm]
\mbox{$SW$-multiplicity map} \mbox{with integer values} (case $b^+>1$)& the
map
$n^\ori_X$ which assigns to a $SC$-structure $W$ the class $n_X(W)$
evaluated with respect to an orientation determined by a choice of a
homology orientation
$\ori$, i.e. an orientation for $H^1(X,\rone)\oplus
H^{2,+}(X,\rone)$.\\[.2cm]
$SW$-structure & a $SC$-structure with non-trivial
$SW$-multiplicity.\\[.2cm]
basic class & the Chern class of a $SW$-structure.
\end{tabular}

\begin{thm}[\cite{b}]
\label{can_inv}
If $X$ is a minimal K\"ahler surface of non-negative Kodaira dimension
then the homology class $k$ of its canonical divisor is invariant up to
sign under oriented diffeomorphism.
\end{thm}

\begin{thm}[\cite{b}]
\label{uniqueSW}
Let $X$ be a K\"ahler surface with positive geometric genus then $k$ and
$-k$ are basic classes. Moreover there is a unique SW-structure with
determinant equal to $-k$ modulo torsion.
\end{thm}

\begin{thm}[\cite{b},\cite{w}]
\label{SWtransf}
The SW-multiplicities with integer values of K\"ahler surfaces with
$b^+\geq3$ have the following properties:
\begin{enumerate}
\item
invariance under diffeomorphisms $\varphi:X\to Y$:
$$
n^\ori_Y(W)=n^{\varphi^\ast\ori}_X(\varphi^\ast W).
$$
\item
antisymmetry under change of homology orientation:
$$
n^\ori_X(W)=-n^{\text{-}\ori}_X(W).
$$
\end{enumerate}
\end{thm}

\begin{lemma}
\label{ori_inv}
Let $\varphi$ be a diffeomorphism of a K\"ahler surface $X$ changing the
homology orientation, $c$ a class fixed by $\varphi$ up to torsion. If $W$
is a SW-structure with determinant equal to $c$ modulo torsion, then so is
$\varphi^\ast W$. Moreover the two structures are different.
\end{lemma}

\proof
Obviously $W$ and $\varphi^\ast W$ have the same determinant up to
torsion. By the preceding theorem and the assumption on $\varphi$ we have:
$$
n^\ori_X(\varphi^\ast W)=-n^{\text{-}\ori}_X(\varphi^\ast W)=
-n^{\varphi^\ast\ori}_X(\varphi^\ast W)=-n^\ori_X(W).
$$
Hence the values of $\sw$ are non-trivial both, and different.
\qed

\begin{lemma}
\label{hyp}
The homology lattice $(\Hbar_2,\pairing)$ of an elliptic surface with
positive Euler number and positive geometric genus contains a hyperbolic
direct summand orthogonal to the homology class $k$ of its canonical
divisor.
\end{lemma}

\proof
This is an obvious corollary to the classification of indefinite even
unimodular forms in the $K3$ case where $k$ is trivial.\\
Otherwise choose a class which pairs to one with a primitive class
commensurable with $k$. Together they span a unimodular sublattice, thus
its orthogonal complement is an orthogonal direct summand of $\Hbar$.
The assumption on the geometric genus implies that the intersection form
restricted to this summand is indefinite. Of course it is unimodular, and
it is even, because $k$ is a characteristic element for the intersection
form. Therefore as in the trivial case we may appeal to the classification
of indefinite even unimodular forms to get a hyperbolic summand as claimed.
\qed

\begin{lemma}
\label{inversion}
In $\Ok$ there is an isometry $\imath$ of the homology lattice $\Hbar$ of
an elliptic surface as above which is not induced from any diffeomorphism
of the surface.
\end{lemma}

\proof
Let $\imath$ be inversion at zero on a hyperbolic summand as given in
lemma \ref{hyp} and the identity on the orthogonal complement, then
$\imath$ stabilizes the canonical class but changes the homology
orientation.\\
Given a diffeomorphism which changes the homology orientation, then
by lemma \ref{ori_inv} a class which is fixed up to torsion does not have a
unique SW-structure. But as we know from thm.\ref{uniqueSW} there is a
unique SW-structure with determinant $-k$ modulo torsion, so we may
reverse the conclusion above to disprove the existence of a diffeomorphism
inducing $\imath$.
\qed
\\[-1cm]

\begin{lemma}[\cite{fm1}]
\label{cpxconj}
The diffeomorphism group of an elliptic surface with positive Euler number
contains a diffeomorphism $\sigma$ with $\sigma(k)=-k$.
\end{lemma}

\noindent
The idea of the construction is to embed the surface in projective space.
Complex conjugation then maps it to another surface, the fibre class to
the negative of the fibre class on the second surface.\\

\noindent
On the other hand the two surfaces are deformation equivalent, a family
containing both gives rise to at least one diffeomorphism between them
preserving the canonical class. Since the canonical class is a multiple of
the fibre class, the composite map does the job.\\[.2cm]

\begin{thm}
Let $X$ be a minimal elliptic surface with positive Euler number, $k$ the
homology class of its canonical divisor and
$$
\psi:\Diff(X)\tto\Orth(\Hbar_2(X))
$$
the natural homomorphism, then the image of $\psi$ is given by
\begin{equation*}
\im\psi=
 \begin{cases}
   \Ospk\adsig & \text{ in case $p_g>0$,}\\
   \quad=\Ospin & \text{ if $X$ is a $K3$ surface},\\[.2cm]
   \:\Ok\adid & \text{ in case $p_g=0$, $X$ non rational,}\\
   \quad=\Orth & \text{ in the special case of Enriques surfaces},\\[.2cm]
   \qquad\Orth & \text{ in case $X$ rational}.\\
 \end{cases}
\end{equation*}
\end{thm}

\proof
The result is well known in the rational case, cf. \cite{wall}, and has
been previously established for $K3$ surfaces by Donaldson \cite{d}.\\
Our proof for $K3$ surfaces is based on the fact
$\Ospin\subset\im\psi\subset\Orth$, which is obvious by thm.\ref{E+K3},
and on the existence of $\imath\in\Orth,\imath\notin\im\psi$ as given by
lemma \ref{inversion}. Then the result follows from $[\Orth\!:\!\Ospin]=2$
(notice that consequently $\sast\in\Ospin$).\\
For Enriques surfaces we have $\Ospin\subset\im\psi\subset\Orth$ as well,
but this time lemma \ref{inversion} does not apply due to the fact that
$p_g=0$, i.e. the canonical divisor - though numerically trivial - is not
effective. Instead we use $\sast\in\im\psi$ of lemma \ref{cpxconj} which
can be shown not to be contained in $\Ospin$. Again by
$[\Orth\!:\!\Ospin]=2$ our claim is proved.\\
In the case of positive geometric genus and $k$ non-trivial we have
$\Ospk\subset\im\psi$ by the corollary to thm.\ref{Xdqm_diffgrp}, since $k$
is a non-trivial multiple of the fibre class.
On the other hand we just got $\im\psi\subset\Ok\adid$ which contains
$\Ospk$ as subgroup of index at most four.
Moreover the preceding lemmas provide us with elements
$\sast,\imath\in\Ok\adid$ such that
$\quad\sast\in\im\psi,\,\sast\notin\Ospk\quad\text{and}
\quad\imath\in\Ok,\,\imath\notin\im\psi.\quad$
Hence $\Ok\adid$ contains $\im\psi=\Ospk\adsig$ as subgroup of index two.\\
The argument is similar in the case of vanishing geometric genus and $X$
non-rational, i.e. with a $ne@!@!f$ canonical divisor, so $k$ is a
non-trivial positive multiple: We have $\Ok\subset\im\psi\subset\Ok\adid$
by the appropriate results, the corollary to thm.\ref{Xdqm_diffgrp} and
thm.\ref{can_inv}. Since $\sast$ from lemma
\ref{cpxconj} is in $\im\psi$ but not in $\Ok$ the second relation is in
fact an equality of the groups.
\qed

\noindent
{\sc Institut f\"ur Mathematik, Universit\"at Hannover,\\
Postfach 6009,
30060 Hannover}\\
{\sl E-mail adress}: loenne@@math.uni-hannover.de


\begin{thebibliography}{F/M1}

\bibitem[BPV]{bpv} Barth, W., Peters, C., Van de Ven, A.: {\sl Compact
Complex Surfaces}, Ergebnisse, 3.Folge, Band 4, Springer, Berlin, etc.
(1984)
\bibitem[B]{b} Brussee, R.: {\sl Some $C^\infty$ properties of K\"ahler
surfaces}, (preprint alg-geom/9503004)
\bibitem[D]{d} Donaldson, S.: {\sl Polynomial invariants for smooth
four-manifolds}, Topology 29(1990), 257-315
\bibitem[E]{e} Ebeling, W.: {\sl The Monodromy Groups of Isolated
Singularities of Complete Intersections}, LNM  1293, Springer, Berlin,
etc. (1987)
\bibitem[E/O]{eo} Ebeling, W., Okonek, C.: {\sl On the diffeomorphism
groups of certain algebraic surfaces}, Enseign. Math. 37(1991), 249-262
\bibitem[F]{f} Freedman, M.: {\sl The topology of four-dimensional
manifolds}, J. Diff. Geom. 17(1982), 357-454
\bibitem[F/M1]{fm1} Friedman, R., Morgan, J.: {\sl On the diffeomorphism
types of certain algebraic surfaces I}, J. Diff. Geom. 27(1988), 297-369
\bibitem[F/M2]{fm2} Friedman, R., Morgan, J.: {\sl Algebraic surfaces and
Seiberg-Witten Invariants}, (preprint alg-geom/9502026)
\bibitem[F/M]{fm} Friedman, R., Morgan, J.: {\sl Smooth Four-Manifolds
and Complex Surfaces}, Ergebnisse, 3.Folge, Band 27, Springer, Berlin, etc.
(1994)
\bibitem[K]{k} Kreck, M.: {\sl Isotopy classes of diffeomorphisms of
$(k-1)$-connected almost parallelizable $2k$-manifolds}, in: Dupont, J.L.,
Madsen, I.H. (ed.) {\sl Algebraic Topology Aarhus 1978}, LNM 763,
Springer, Berlin, etc. (1979), 643-663
\bibitem[La]{la} Lamotke, K.: {\sl The topology of complex projective
varieties after S. Lefschetz}, Topology 20(1981), 15-51
\bibitem[Q]{q} Quinn, F.: {\sl Isotopy of 4-manifolds}, J. Diff. Geom.
24(1986), 343-372
\bibitem[Wa]{wall} Wall, C.T.C.: {\sl Diffeomorphisms of 4-manifolds}, J.
LMS 39(1964), 131-140
\bibitem[W]{w} Witten, E.: {\sl Monopoles and Four-Manifolds}, (preprint
hep-th/9411102), Math. Res. Letters 1(1994), 769-796

\end{thebibliography}
\end{document}